\documentclass[final,a4paper,leqno]{siamltex}

\usepackage[english]{babel}
\usepackage{amsmath,amssymb,amsfonts,graphicx,subcaption,hyperref}
\usepackage[ruled]{algorithm2e} 
\usepackage{myshortcuts}
\usepackage{xr}

\usepackage{pgfplots}
\usetikzlibrary{external}
\usetikzlibrary{decorations.pathreplacing}

\iffalse 
\usepackage{tikz}
\usepackage{tikzscale}
\pgfplotsset{
    compat=newest,
    tick label style={font=\scriptsize},
    label style={font=\scriptsize},
    legend style={font=\scriptsize}
}
\usetikzlibrary{decorations.pathreplacing}
\iftrue 
\usepgfplotslibrary{external} 
\else
\renewcommand{\tikzsetnextfilename}[1]{}
\fi
\else
\usepackage{tikzexternal}
\tikzsetexternalprefix{fig/}

\fi

\title{Linearizability of eigenvector nonlinearities \thanks{The work by Rob Claes and Karl Meerbergen is supported by the Research Foundation – Flanders (FWO) Grant G0B7818N and the KU Leuven Research Council.}}

\author{Rob Claes\thanks{Department of Computer Science, KU Leuven, Leuven, Belgium, email: \texttt{\{rob.claes, karl.meerbergen\}@kuleuven.be}} \and Elias Jarlebring \thanks{Department of Mathematics, Royal Institute of Technology (KTH), Stockholm, SeRC Swedish e-Science Research Center, email: \texttt{\{eliasj, pup\}@kth.se}} \and Karl Meerbergen \footnotemark[2] \and Parikshit Upadhyaya \footnotemark[3]}

\date{\today}
\begin{document}

\maketitle
\begin{abstract}
    We present a method to linearize, without approximation, a specific class of eigenvalue problems with eigenvector nonlinearities (NEPv), where the nonlinearities are expressed by scalar functions that are defined by a quotient of linear functions of the eigenvector.
    The exact linearization relies on an equivalent multiparameter problem (MEP) that contains the exact solutions of the NEPv. Due to the characterization of MEPs in terms of a generalized eigenvalue problem this provides a direct way to compute all NEPv solutions for small problems, and it opens up the possibility to develop locally convergent iterative methods for larger problems.
    Moreover, the linear formulation allows us to easily determine the number of solutions of the NEPv.
    We propose two numerical schemes that exploit the structure of the linearization: inverse iteration and residual inverse iteration.
    We show how symmetry in the MEP can be used to improve reliability and reduce computational cost of both methods.
    Two numerical examples verify the theoretical results and a third example shows the potential of a hybrid scheme that is based on a combination of the two proposed methods.
\end{abstract}

\begin{keywords}
    nonlinear eigenvalue problem, multiparameter eigenvalue problem, linearization
\end{keywords}

\begin{AMS}
    65F15, 65H17
\end{AMS}

\section{Introduction}\label{sec:intro}

Consider the following eigenvalue problem
\begin{equation}  \label{eq:nepv_general}
(A+\lambda B+f_1(x)C_1+\cdots+f_m(x)C_m)x=0
\end{equation}
with $A,B\in \mathbb{C}^{n\times n}$ and eigenvector nonlinearities defined by $f_i: \mathbb{C}^n\to \mathbb{C}$ and $C_i\in\mathbb{C}^{n\times n}$ for $i=1,2,\dots,m$.
The goal is to find $\lambda\in\mathbb{C}$ and $x\in\mathbb{C}^n\backslash\{0\}$ such that the eigenpair $(\lambda,x)$ solves \eqref{eq:nepv_general}.
This problem belongs to the class of nonlinear eigenvalue problems with eigenvector nonlinearities or \emph{NEPv}.
The functions $f_i$ are defined such that $f_i(x)=f_i(\alpha x)$ for any $\alpha\in\mathbb{C}\backslash\{0\}$.

Problems with eigenvector nonlinearities arise in a wide range of fields. The most common application is in the field of quantum chemistry for electronic structure calculations \cite{Szabo:1996:QC}, where numerical discretizations of the Schr\"{o}dinger equation lead to a NEPv. See \cite{Saad2010} for a summary of numerical methods and discretization approaches. In these problems, the nonlinearity may depend on a set of eigenvectors, instead of a single eigenvector as in \eqref{eq:nepv_general}. The Gross-Pitaevskii equation is a related application in quantum physics, where the goal is to compute the ground state of bosons \cite{Bao2004, Jarlebring2014Inverse}. A widely used approach for the solution in both these fields is the self-consistent field (SCF) iteration \cite{Cances:2000:SCF, Liu2015} and its variants \cite{Pulay:1980:CONVERGENCE, Rohwedder2011}. SCF is an iterative scheme that involves solving a linear eigenproblem at each step using the iterate from the previous step, until convergence. For convergence results on SCF, see \cite{Yang:2009:SCF, Levitt2012, Liu2014, Upadhyaya2019}.

A different class of applications is in the field of data science. Spectral clustering using the eigenpairs of the p-Laplacian relies on solving a NEPv \cite{Hein:2009:PLAPLACIAN}. In \cite{Hein:2010:IPM}, the authors propose an inverse power method for the problem. A Rayleigh quotient minimization approach applied to Fisher linear discriminant analysis is discussed in \cite{Bai_LV:2018}. In \cite{Tudisco:2019:CORE}, an equivalence between a NEPv and a new model for core-periphery detection in network science is shown.

Linearization without approximation
in the sense of companion linearizations is a crucial
tool for eigenvalue nonlinearities \cite{mackey2006structured}.
It is not only used as a procedure in itself,
but serves as a basis of projection
methods, where the linearization can provide insight in
building a projection subspace. In this paper, we will study
linearization in this sense, and for simplicity
the term \emph{linearization} will be used to refer to this, and not
to the linear approximation of a nonlinear problem.

In this paper, we extend this notion of linearization to eigenvector nonlinearities by focusing on a subset of problems that can be written as \eqref{eq:nepv_general} with
\begin{equation} \label{eq:nepv_general_f}
f_i(x) = \frac{r_i^Tx}{s_i^Tx}, \quad i=1,2,\dots,m
\end{equation}
and $r_i,s_i\in \mathbb{C}^n$.
To illustrate the concept of the new method, consider the following example for $m=1$:
\begin{gather*}\label{eq:MEP2}
(A+\lambda B +f_1(x)C_1)x = 0 \\
f_1(x) = \frac{r_1^Tx}{s_1^Tx}.
\end{gather*}
The linearization then consists of writing this problem as a linear multiparameter eigenvalue problem (MEP) as follows:
\begin{equation}\label{eq:intro_ex}
\begin{cases}
-Ax = (\lambda B+\mu C_1)x \\
-(A+g_1r_1^T)x = (\lambda B + \mu(C_1-g_1s_1^T))x
\end{cases}
\end{equation}
with mild conditions on $g_1\in \mathbb{C}^n\backslash\{0\}$.
Subtracting the first equation from the second results in
\begin{equation*}
-g_1r_1^Tx = -\mu g_1s_1^Tx,
\end{equation*}
and since $g_1\ne 0$ it follows that
\begin{equation*}
\mu = \frac{r_1^Tx}{s_1^Tx} = f_1(x).
\end{equation*}
Finally, substituting this value of $\mu$ in the first equation of \eqref{eq:intro_ex} shows that solutions of the NEPv are contained in the linear MEP without approximation.

In this paper we provide theoretical results concerning this
class of problems, as well as new algorithms. More specifically,
we describe the procedure outlined above for the general case,
including a charactarization of a linearization to an MEP,
which subsequently can be linearized to a generalized linear
eigenvalue problem (GEP) using standard MEP techniques.
We first show that, in general, the number of solutions to a problem belonging to this class is bounded.
Theory is provided that justifies the method, in particular that
the MEP is equivalent to the NEPv under certain conditions. The
MEP always contains the NEPv solutions, but depending on certain
free parameters (later denoted $g_1,\ldots,g_m$),
the converse is not always true. We derive sufficient conditions which are generically satisfied,
and show examples that illustrate degenerate cases.

The theory also reveals that it is sufficient to
consider MEP solutions with a particular structure,
namely that the eigenvector components are identical (referred to as symmetric MEP solutions).
It turns out that certain iterative MEP methods preserve the symmetry,
if they are initiated with a symmetric initial vector.
This allows us to improve the algorithms. In particular,
we improve residual inverse iteration \cite{plestenjak2016numerical} and inverse iteration applied to the GEP. We provide
convergence theory illustrating how the convergence theory
for the general methods can be specialized.
A set of experiments highlights the strengths of each method,
and shows how the structure can arise from a PDE-problem.

The remainder of this paper is organized as follows.
Section~\ref{sec:prelim} gives an introduction to the theory of MEPs and shows how operator determinants can be used to extract the eigenpairs.
This theory is used in Section~\ref{sec:lin_theory} to demonstrate the new linearization method that defines a multiparameter eigenvalue problem containing solutions to the nonlinear problem.
Section~\ref{sec:num} is devoted to numerical methods, convergence
characterizations and simulations.

\section{Preliminaries} \label{sec:prelim}

The proposed method transforms the NEPv to a linear multiparameter eigenvalue problem.
This section therefore first summarizes the relevant theory regarding multiparameter eigenvalue problems and their relation to operator determinants based on \cite{Atkinson:1972:MULTIPARAMETER, Kosir:1994:FINITE}.

We write the multiparameter eigenvalue problem (MEP) as a set of $l$ coupled equations
\begin{equation} \label{eq:MEP}
\begin{cases}
V_{1,0}x_1 &= (\mu_1 V_{1,1}+\mu_2 V_{1,2} + \cdots + \mu_{l} V_{1,l})x_1 \\
\quad\vdots \\
V_{l,0}x_{l} &= (\mu_1 V_{l,1}+\mu_2 V_{l,2} + \cdots + \mu_{l} V_{l,l})x_{l}
\end{cases}
\end{equation}
with $V_{i,j}\in \mathbb{C}^{n_i\times n_i}$, $\mu_i\in\mathbb{C}$ and $x_i\in\mathbb{C}^{n_i}$ for $i=1,2,\dots,l$ and $j=0,1,\dots,l$.
The goal is to find tuples $(\mu_1,\dots,\mu_{l},x_1,\dots,x_{l})$ with all $x_i$ nonzero that solve \eqref{eq:MEP}.
In general, this problem has $\prod_{i=1}^{l}n_i$ different isolated solution tuples.
An important theoretical result is that these tuples can be extracted from a linear generalized eigenvalue problem (GEP).
This equivalence between the MEP and the GEP is described by operator determinants.
The remainder of this section describes such operator determinants and shows how they are used to construct an equivalent GEP to a given MEP.

In general, an operator determinant of a $k\times k$-array of matrices $D_{i,j}$ is denoted by
\begin{equation}\label{eq:det}
    \begin{vmatrix}D_{1,1}&\cdots &D_{1,k}\\
    \vdots &\ddots&\vdots\\
    D_{k,1}&\cdots &D_{k,k}
    \end{vmatrix}
    _\kron
\end{equation}
where $D_{i,j}\in\CC^{n_i\times n_i}$, $i,j=1,\dotsc,k$.
The $\kron$ subscript corresponds to usage of
standard formulas for determinants, but  where
multiplication is replaced by the Kronecker product
and the order of operations is chosen consistently.
The operator determinant can therefore be defined by
an operator cofactor expansion, e.g.,
expansion by the first block row yields
\begin{equation}
   \begin{aligned}
   \begin{vmatrix}D_{1,1}&\cdots &D_{1,k}\\
   \vdots &\ddots&\vdots\\
   D_{k,1}&\cdots &D_{k,k}
   \end{vmatrix}
   _\kron:=&
   D_{1,1}\otimes
   \begin{vmatrix}D_{2,2}&\cdots &D_{2,k}\\
   \vdots &\ddots&\vdots\\
   D_{k,2}&\cdots &D_{k,k}
   \end{vmatrix}
   _\kron +\ldots \\
   &+(-1)^{k+1}
   D_{1,k}\otimes
   \begin{vmatrix}D_{2,1}&\cdots &D_{2,k-1}\\
   \vdots &\ddots&\vdots\\
   D_{k,1}&\cdots &D_{k,k-1}
   \end{vmatrix}
   _\kron.
   \end{aligned}
\end{equation}

Notice that not all properties for standard determinants can be applied to operator determinants.
The following theorem lists properties that hold for operator determinants.
\begin{theorem}
    \cite[Theorem 6.2.1]{Atkinson:1972:MULTIPARAMETER}
    For a determinant \eqref{eq:det}, we have the following properties
    \begin{enumerate}
        \item If two columns are interchanged, the sign of the determinant is reversed.
        \item If two columns are identical, the determinant vanishes.
        \item The value of the determinant is unchanged if a scalar multiple of one column is added to another column.
    \end{enumerate}
\end{theorem}
An important difference to standard determinants, is that in general these properties only hold for columns and not for rows.

Now define $\Delta_0$ and $\Delta_i$, $i=1,\dots,l$ as operator determinants constructed from the matrices of the MEP \eqref{eq:MEP}:
\begin{align}  \label{eq:Delta01}
\Delta_0&=
\begin{vmatrix}V_{1,1}&V_{1,2}&\cdots &V_{1,l}\\
\vdots &\vdots & &\vdots\\
V_{l,1}&V_{l,2}&\cdots &V_{l,l}
\end{vmatrix}
_\kron,\\
\Delta_i&=
\begin{vmatrix}V_{1,1} & \cdots & V_{1,i-1} & V_{1,0} & V_{1,i+1} & \cdots & V_{1,m+1}\\
\vdots &&\vdots &\vdots &\vdots & & \vdots\\
V_{l,1} & \cdots & V_{l,i-1} & V_{l,0} & V_{l,i+1} & \cdots & V_{l,l}
\end{vmatrix}_\otimes.
\end{align}
We say that the MEP \eqref{eq:MEP} is nonsingular if and only if $\Delta_0$ is nonsingular.
The following theorem summarizes the relation between the eigenvalues of the MEP and its operator determinants.
\begin{theorem}\cite[Chapter~6]{Atkinson:1972:MULTIPARAMETER}\label{thm:gen_eig}
    The multiparameter eigenvalue problem \eqref{eq:MEP} is equivalent to the following system of generalized eigenvalue problems
    \begin{equation}\label{eq:op_det1}
    \Delta_iz = \mu_i\Delta_0z, \quad i=1,\dots,l\\
    \end{equation}
    with $z=x_1\otimes\dots\otimes x_{l}$.
    Moreover if the MEP is nonsingular, all matrices $\Delta_0^{-1}\Delta_i$ commute such that the vector $z$ always exists.
\end{theorem}

This theorem suggests a way to extract the eigenvalues and eigenvectors of this system.
In short, it suffices to solve one eigenvalue problem, e.g. $(\Delta_1,\Delta_0)$, to find the complete eigentuples $(\mu_1,\ldots,\mu_l,x_1,\ldots,x_l)$.
It is possible to reduce the computational cost by exploiting the structure of the operator determinants \cite{meerbergen2015sylvester} or to use a subspace method for larger problems \cite{Hochstenbach:2005:TWOPARAMETER, Hochstenbach:2002:JDDEFTWOPARAM,hochstenbach2019subspace}.
Section \ref{sec:lin_nepv_gep} discusses the system \eqref{eq:op_det1} in greater detail.

The linearization approach presented in this paper heavily relies on the concept of symmetric decomposable vectors.
\begin{definition}\label{def:symdec}
    A symmetric decomposable vector of rank $k$ is a vector $z\in\mathbb{C}^{n^l}$ that can be decomposed in $k$ linearly independent vectors as $z = \sum_{j=1}^{k} \alpha_jx_j\otimes\dots\otimes x_j$ with $x_j\in\mathbb{C}^n$.
\end{definition}

For the remainder of this paper, a symmetric eigenvector is thus defined as a symmetric decomposable vector of rank 1 that is a solution to \eqref{eq:op_det1}, and $\lambda$ is its corresponding symmetric eigenvalue.

We also use the term \emph{symmetric} in the context of solutions to MEPs.
\begin{definition}\label{def:simpeig}
A solution $(\mu_1,\ldots,\mu_l,x_1,\ldots,x_l)$ of the MEP \eqref{eq:MEP} is called symmetric if
  \[
 x_i = \alpha_i x_1,\quad  \forall i\in\{2,\ldots,l\},\quad \alpha_i \in\mathbb{C}\setminus\{0\}.
  \]
\end{definition}
If the assumptions of Theorem~\ref{thm:gen_eig} are satisfied, then a symmetric solution of the MEP \eqref{eq:MEP}, for example, $(\mu_1,\ldots,\mu_l,x_1,\ldots,x_1)$ corresponds to an a rank one symmetric eigenvector $z = x_1\kron x_1\kron \cdots \kron x_1$ of the equations in \eqref{eq:op_det1}.

\section{Linearization of NEPv} \label{sec:lin_theory}
Before introducing the linearization of the NEPv \eqref{eq:nepv_general}-\eqref{eq:nepv_general_f}, we derive the number of solutions of this problem.
\begin{theorem}\label{thrm:nb_solutions}
    The NEPv \eqref{eq:nepv_general}-\eqref{eq:nepv_general_f} has at most $N_s=\begin{pmatrix}n+m\\m+1\end{pmatrix}$ isolated solutions.
\end{theorem}
\begin{proof}
    First, rewrite the NEPv \eqref{eq:nepv_general}-\eqref{eq:nepv_general_f} as
    \begin{subequations} \label{eq:nb_solutions_system}
        \begin{align}
        (\lambda_0A + \lambda_1B + \lambda_2C_1 + \cdots + \lambda_{m+1}C_{m+1})x = 0 & & \\
        \lambda_0r_i^Tx - \lambda_{i+1}s_i^Tx = 0 && i=1,\dots,m
        \end{align}
    \end{subequations}
    with $\lambda_k\in\mathbb{C}$ for $k=0,\dots,m+1$.
    This is a system of $n+m$ bihomogenous equations of degree $(1,1)$ on the multiplicative projective space $\mathbb{P}^{m+1}\times\mathbb{P}^{n-1}$.
    Here, a projective space $\mathbb{P}^k$ is defined by complex coordinates $(Z_1,\ldots,Z_{k+1})\in\mathbb{C}^{k+1}\setminus\{0\}$ such that $(Z_1,\ldots,Z_{k+1})$ and $(\gamma Z_1,\ldots,\gamma Z_{k+1})$ correspond to the same point in the projective space for every $\gamma\in\mathbb{C}\setminus\{0\}$.
    By the multi-homogenous B\'ezout theorem \cite{shafarevich2013basic}, the number of solutions is bounded from above by the coefficient of $t_1^{m+1}t_2^{n-1}$ in the expansion of $(t_1+t_2)^{n+m}$, which is the binomial coefficient $\begin{pmatrix}n+m\\m+1\end{pmatrix}$.
    The finite solutions of the NEPv correspond to the solutions of \eqref{eq:nb_solutions_system} by the equalities $\frac{\lambda_{i+1}}{\lambda_0} = \frac{r_i^Tx}{s_i^Tx} = f_i(x)$ for $i=1,\dots,m$ if  $\lambda_0\neq0$ and $s_i^Tx\neq0$, which means that \eqref{eq:nepv_general}-\eqref{eq:nepv_general_f} cannot have more than $N_s=\begin{pmatrix}n+m\\m+1\end{pmatrix}$ finite solutions.
\end{proof}

It is important to emphasize that for generic problems of type \eqref{eq:nepv_general}-\eqref{eq:nepv_general_f}, the upper bound $N_s$ is almost always reached.
For general coefficients $A,B,C_i,r_i,s_i$ for $i=1,\dots,m$, the system \eqref{eq:nb_solutions_system} has exactly $N_s$ solutions with probability 1.
It also holds generically that $\lambda_0\neq0$ and $s_i^Tx\neq0$, which means that all $N_s$ solutions of \eqref{eq:nb_solutions_system} correspond to finite isolated solutions of \eqref{eq:nepv_general}-\eqref{eq:nepv_general_f} and the upper bound is reached.

\subsection{Linearization of NEPv to MEP}
We will show that the original problem NEPv \eqref{eq:nepv_general}
can be cast to an equivalent MEP.
We say that we linearize the NEPv to this MEP
\begin{subequations}\label{eq:nepv_general_linearized}
    \begin{eqnarray}    -Ax_1 &=& (\lambda B + \mu_1C_1 + \cdots + \mu_mC_m)x_1 \label{eq:nepv_general_linearized_a}\\
    -(A+g_1r_1^T)x_2 &=& (\lambda B + \mu_1(C_1-g_1s_1^T) + \cdots + \mu_mC_m)x_2 \label{eq:nepv_general_linearized_b}\\
    &\vdots& \notag\\
    -(A+g_mr_m^T)x_{m+1} &=& (\lambda B + \mu_1C_1 + \cdots + \mu_m(C_m-g_ms_m^T))x_{m+1}
    \end{eqnarray}
\end{subequations}
where $g_1,\ldots,g_m\neq 0$ can be arbitrary vectors under mild conditions that are further discussed in Section~\ref{subsec:choiceg}.

\begin{theorem}\label{thm:nepvtomep}
  Assume $(\lambda,x)$ is a solution to \eqref{eq:nepv_general}, then
\[
(\lambda, \mu_1, \cdots, \mu_m, x_1,\ldots,x_{m+1})=
(\lambda, f_1(x), \cdots, f_m(x), x,\ldots,x)
\]
is a solution to \eqref{eq:nepv_general_linearized}.\\
\end{theorem}
\begin{proof}
    Define the variables $x_1,\ldots,x_{m+1}$ and
    $\mu_1,\ldots,\mu_m$ from the solution to \eqref{eq:nepv_general}
    as follows:
    $\mu_1:=f_1(x)$,
    $\ldots$,
    $\mu_m:=f_m(x)$
    and $x_1=\ldots=x_{m+1}=x$.
    Equation \eqref{eq:nepv_general_linearized_a} is
    satisfied due to \eqref{eq:nepv_general}.
    Equation \eqref{eq:nepv_general_linearized_b}
    is satisfied by forming the sum of \eqref{eq:nepv_general}
    and the following identity multiplied by $g_1$
    \[
    r_1^Tx=\mu_1s_1^Tx,
    \]
    which follows from the definition of $\mu_1=f_1(x)=\frac{r_1^Tx}{s_1^Tx}$.
    The other equations in \eqref{eq:nepv_general_linearized}
    follow from an analogous construction.
\end{proof}

 The converse of the theorem above is illustrated by two theorems that follow. The first theorem illustrates that every symmetric solution to the MEP leads to a solution to the NEPv under an assumption on $x$, which characterizes potential spurious solutions introduced by our approach.

\begin{theorem}\label{thm:equiv2}
    Assume $(\lambda,\mu_1,\ldots,\mu_m, x_1,\cdots, x_{m+1})$ is a solution to \eqref{eq:MEP} where $V_{ij}$ are given by \eqref{eq:nepv_general_linearized}.
    If $x_1=x_2=\cdots=x_{m+1}=x$ then $(\lambda, x_1)=(\lambda,x)$ is
    a solution to \eqref{eq:nepv_general},
    if $s_j^Tx\neq 0$, for all $j=1,\ldots,m$.
\end{theorem}
\begin{proof}
    By subtracting equation
    \eqref{eq:nepv_general_linearized_b}
    from
    \eqref{eq:nepv_general_linearized_a} we find that
    \[
    g_1r_1^Tx=\mu_1g_1s_1^Tx.
    \]
    By assumption $s_1^Tx\neq 0$ and $g_1\neq 0$, such that $\mu_1=\frac{r_1^Tx}{s_1^Tx}=f_1(x)$.
    By considering differences between the first and $k$th
    row in \eqref{eq:nepv_general_linearized}
    we obtain analogously
    \begin{equation}  \label{eq:proof_mudef}
    \mu_i=f_i(x),\;\;i=1,\ldots,m.
    \end{equation}
    The connection with \eqref{eq:nepv_general}
    is established by using \eqref{eq:nepv_general_linearized}
    and the formulas for $\mu$ in \eqref{eq:proof_mudef}.
\end{proof}

For the next result, first consider the following problem, with the same eigenvalues as \eqref{eq:nepv_general_linearized}.
\begin{subequations}\label{eq:nepv_general_linearizedt}
    \begin{eqnarray}
    -A^Ty_1 &=& (\lambda B^T + \mu_1C_1^T + \cdots + \mu_mC_m^T)y_1 \\
    -(A^T+r_1g_1^T)y_2 &=& (\lambda B^T + \mu_1(C_1^T-s_1g_1^T) + \cdots + \mu_mC_m^T)y_2 \\
    &\vdots& \notag\\
    -(A^T+r_mg_m^T)y_{m+1} &=& (\lambda B^T + \mu_1C_1^T + \cdots + \mu_m(C_m^T-s_mg_m^T))y_{m+1}
    \end{eqnarray}
\end{subequations}
A tuple of vectors $(y_1, \ldots, y_{m+1})$ that solves this problem is called a left eigenvector of \eqref{eq:nepv_general_linearized}.
The next theorem shows that even non-symmetric solutions to the MEP lead to a solution to the NEPv under an additional assumption on the $g_i$ vectors.
\begin{theorem}\label{thm:always_sym}
    Assume $(\lambda,\mu_1,\ldots,\mu_m, x_1,\cdots, x_{m+1})$ is a solution to \eqref{eq:MEP} where $V_{ij}$ is given by \eqref{eq:nepv_general_linearized}. Let $(\lambda,\mu_1,\ldots,\mu_m,y_1,\ldots,y_{m+1})$ be the solution to the problem given by \eqref{eq:nepv_general_linearizedt}.  If $s_j^Tx_1\neq 0$ and $g_j^Ty_{j+1} \neq 0$ for all $j=1,\ldots,m$, then $(\lambda, x_1)$ is a solution to \eqref{eq:nepv_general},
\end{theorem}
\begin{proof}
Left-multiplying all the subequations except the first in \eqref{eq:nepv_general_linearizedt} by $x_1^T$ and using \eqref{eq:nepv_general_linearized_a}, we have
\begin{equation}
x_1^Tr_jg_j^Ty_{j+1} = \mu_{j}x_1^Ts_jg_j^Ty_{j+1},\quad j=1,\ldots,m.
\end{equation}
The assumptions $s_j^Tx_1 \neq 0$ and $g_j^Ty_{j+1} \neq 0$ for all $j=1,\ldots,m$ lead to
\begin{equation}
\mu_j = \frac{r_j^Tx_1}{s_j^Tx_1},\quad j=1,\ldots,m.
\end{equation}
Replacing the above for $\mu_j$ in the first subequation of \eqref{eq:nepv_general_linearized}, we see that $(\lambda,x_1)$ is a solution to \eqref{eq:nepv_general}.
\end{proof}

Theorem~\ref{thm:always_sym} together with Theorem~\ref{thm:nepvtomep} leads to the conclusion that under certain conditions, for every non-symmetric solution $(\lambda,\mu_1,\ldots,\mu_m,x_1\ldots,x_m)$ to the MEP \eqref{eq:MEP} with the $V_{ij}$ given by \eqref{eq:nepv_general_linearized}, there is a corresponding symmetric solution $(\lambda,\mu_1,\ldots,\mu_m,x_1,\ldots,x_1)$. Since our main interest is to compute solutions to the NEPv, it suffices to restrict our search to the set of symmetric solutions of the MEP. This observation forms the basis of the algorithms in Section~\ref{sec:num}.

\subsection{Linearization of NEPv to GEP} \label{sec:lin_nepv_gep}

\newlength{\bracewidth}
\newcommand{\myunderbrace}[2]{\settowidth{\bracewidth}{$#1$}#1\hspace*{-1\bracewidth}\smash{\underbrace{\makebox{\phantom{$#1$}}}_{#2}}}
The theorems of the preceding subsection show how a solution to the NEPv corresponds to a solution of MEP \eqref{eq:nepv_general_linearized}, and vice versa. Theorem~\ref{thm:gen_eig} shows an equivalence between a MEP and a system of generalized eigenvalue problems (GEP) involving operator determinants, under certain assumptions. For completeness, we now specify the structure of these operator determinants arising from the MEP in \eqref{eq:nepv_general_linearized}. This means that the NEPv can be linearized to a system of GEPs.

The operator determinants of \eqref{eq:op_det1} for the MEP in \eqref{eq:nepv_general_linearized} are as follows:
\begin{equation}
\Delta_0 = \begin{vmatrix}B& C_1& \ldots &C_m \\
            B& C_1-g_1s_1^T& \ldots& C_m\\
            \vdots& \vdots& \ddots& \vdots\\
            B& C_1& \ldots& C_m-g_ms_m^T,
            \end{vmatrix}_\otimes
\end{equation}
and
\begin{equation}
    \Delta_i = \begin{vmatrix}B& C_1& \ldots & -A & \ldots &C_m \\
    B& C_1-g_1s_1^T& \ldots& -A-g_1r_1^T & \ldots & C_m\\
    \vdots& \vdots& \ddots& \vdots &\ddots & \vdots\\
    B& C_1& \ldots& \myunderbrace{-A-g_mr_m^T}{\text{column } i} & \ldots & C_m-g_ms_m^T
    \end{vmatrix}_\otimes
    \vspace{10pt}
\end{equation}
for $i=1,\ldots,m+1$. Note that $\Delta_i$ is obtained by taking $\Delta_0$ and only changing the $i$-th column.

For the linearized NEPv, or the MEP of \eqref{eq:nepv_general_linearized}, the equations of \eqref{eq:op_det1} lead to an equivalent system of $m+1$ GEPs.
\begin{subequations}\label{eq:delta_nepv}
\begin{eqnarray}
\Delta_1 z &=& \lambda \Delta_0 z, \\
\Delta_{i+1} z &=& \mu_{i}\Delta_0 z,\quad i=1,\ldots,m
\end{eqnarray}
\end{subequations}
As discussed at the end of Section~\ref{sec:prelim}, the set of equations in \eqref{eq:delta_nepv} can be solved by only solving in two steps. 
Firstly, the GEP in the first subequation can be solved to find $\lambda$ and $z$. The other values of the eigenvalue tuple can be obtained using the corresponding Rayleigh quotients as
\begin{equation}
\mu_{i} = \frac{z^H \Delta_{i+1}z}{z^H\Delta_0z},\quad i=1,\ldots,m.
\end{equation}
If $z$ is decomposable as $x_1\otimes \ldots \otimes x_{m+1}$ and the assumptions of Theorem~\ref{thm:always_sym} hold, then $(\lambda,x)$ is a solution to the NEPv and we can compute the values of $\mu_i$ without using Rayleigh quotients, as follows:
\begin{equation}
\mu_{i} = \frac{r_{i}^Tx_1}{s_{i}^Tx_1},\quad i=1,\ldots,m.
\end{equation}

We illustrate this with a small example with matrices from $\mathbb{R}^{2\times 2}$, where $m = 1$. Let
\begin{equation}
A = \begin{bmatrix}1& 1\\0& 1\end{bmatrix},\quad B = \begin{bmatrix}1& 2\\3& 4\end{bmatrix},\quad C = \begin{bmatrix}2& 0\\0& 1\end{bmatrix},\quad r = \begin{bmatrix}3\\ 2\end{bmatrix},\quad s = \begin{bmatrix}4\\ 3\end{bmatrix}
\end{equation}
and consider the problem of finding $(\lambda,x)$ in $\mathbb{C}\times\mathbb{C}^2$ such that
\begin{equation}\label{eq:small_ex}
\left(A+\lambda B + \frac{r^Tx}{s^Tx}C\right)x = 0.
\end{equation}
We use $g = \begin{bmatrix}1& 3\end{bmatrix}^T$, which leads to
\begin{subequations}
\begin{eqnarray}
\Delta_ 0 &=& \begin{vmatrix}B& C\\B& (C-gs^T)\end{vmatrix}_{\otimes} = \begin{bmatrix}-4& -7& -4& -6\\ -18& -16& -24& -16\\-6& -9& -9& -14\\-36& -24& -51& -36\end{bmatrix},\nonumber\\
\Delta_1 &=& \begin{vmatrix}-A& C\\-(A+gr^T)& C-gs^T\end{vmatrix}_{\otimes} = \begin{bmatrix}10& 9& 2&3\\30& 22& 12& 8\\0& 0& 6& 6\\0& 0& 21& 15\end{bmatrix}.\nonumber
\end{eqnarray}
\end{subequations}
The four eigenpairs obtained by solving the first GEP in \eqref{eq:delta_nepv} are
\begin{subequations}
\begin{eqnarray*}
\lambda_1 &\approx 5.2462, \quad z_1  \approx \begin{bmatrix}-0.8232\\0.5677\end{bmatrix}\kron \begin{bmatrix}-0.8232\\0.5677\end{bmatrix},\\
\lambda_2 &\approx -0.4224, \quad z_2 \approx \begin{bmatrix}-0.5637\\ 0.8260\end{bmatrix}\kron \begin{bmatrix}-0.5637\\ 0.8260\end{bmatrix},\\
\lambda_3 &\approx -0.4367, \quad z_3 \approx \begin{bmatrix}-0.0672\\ 0.9977\end{bmatrix}\kron \begin{bmatrix}-0.0672\\ 0.9977\end{bmatrix},\\
\lambda_4 &\approx -1.2500, \quad z_4 \approx \begin{bmatrix} -0.4706\\0.8824\end{bmatrix}\kron \begin{bmatrix}0.7682\\ -0.6402\end{bmatrix}.
\end{eqnarray*}
\end{subequations}
It can be verified by substituting into \eqref{eq:nepv_general} that the first three eigenpairs with symmetric eigenvectors give us solutions to the original NEPv.

However, the fourth eigenpair is a spurious solution in the sense that it does not lead to a solution of the original NEPv, even though it is a solution to the MEP. The fourth eigenvalue can be analytically verified to be $\lambda = -\frac{5}{4}$ with left eigenvector
\begin{equation}
w_4 = \begin{bmatrix}-1\\ \frac{1}{3} \end{bmatrix}\kron \begin{bmatrix}-1\\ \frac{1}{3}\end{bmatrix}.
\end{equation}
Hence, 
left eigenvectors of the MEP are given by
$y_1 = y_2 = \left(-1, \frac{1}{3}\right)^T$. This implies that $g^Ty_1 = g^Ty_2 = 0$, which is a violation of the assumptions of Theorem~\ref{thm:always_sym}, which would otherwise state that the $-5/4$ is an eigenvalue of the NEPv.
This explains why the fourth eigenpair corresponds to a spurious solution. These observations are confirmed by the following two lemmas that will help to understand the numerical methods.

\begin{lemma} \label{lmm:rankone_symmetric}
    Under the assumptions of Theorem~\ref{thm:gen_eig} and Theorem~\ref{thm:always_sym}, eigenvectors of simple eigenvalues of \eqref{eq:delta_nepv} are rank one symmetric.
\end{lemma}
\begin{proof}
    Consider a simple eigenpair $(\lambda,z_1)$ of \eqref{eq:delta_nepv} such that that $z_1 = x_1\kron x_2\kron\ldots\kron x_{m+1}$. If the assumptions of Theorem~\ref{thm:gen_eig} and Theorem~\ref{thm:always_sym} hold, then $(\lambda,x_1)$ is a solution to \eqref{eq:nepv_general_linearized} and hence, $\lambda$ is not a spurious eigenvalue. Then, Theorem~\ref{thm:nepvtomep}, dictates that $(\lambda,f_1(x_1),\ldots,f_m(x_1),x_1,\ldots,x_1)$ is a solution to \eqref{eq:nepv_general_linearized}. Therefore, $(\lambda,z_2)$ is also an eigenpair where $z_2 = x_1\kron x_1\kron\cdots \kron x_1$. Since $\lambda$ is a simple eigenvalue, $z_1$ and $z_2$ are linearly dependent, which implies
    \begin{equation*}
    x_i = \alpha_i x_1,\quad i=2,\ldots,m+1
    \end{equation*}
    for $\alpha_i \in\mathbb{C}$. Hence, $z_1$ is rank one symmetric. 
\end{proof}

\begin{lemma} \label{lmm:symmetric_vectors}
    If the NEPv \eqref{eq:nepv_general}-\eqref{eq:nepv_general_f} has $N_s$ distinct eigenvalues, the operator determinant system \eqref{eq:delta_nepv} has $N_s$ linearly independent rank one symmetric eigenvectors.
\end{lemma}
\begin{proof}
    If \eqref{eq:nepv_general}-\eqref{eq:nepv_general_f} has $N_s$ distinct eigenvalues, the operator determinant system \eqref{eq:delta_nepv} has $N_s$ symmetric eigenvalues, and $N_s$ corresponding rank one symmetric eigenvectors.
    Since these eigenvectors belong to distinct eigenvalues of \eqref{eq:delta_nepv}, they are linearly independent.
\end{proof}

Note that the space of symmetric vectors has dimensions $N_s$, which means that there can be no more than $N_s$ linearly independent rank one symmetric vectors.

\subsection{Remarks on the choice of $g_i$}\label{subsec:choiceg}
For a NEPv with $m$ different nonlinear terms $f_i$, our linearization technique introduces $m$ additional vectors $g_i$ that can be chosen freely. The theorems in the previous section suggest that almost all choices of $g_i$ are effective, in the sense that resulting MEP can be used to find solutions to the original NEPv. More specifically, Theorem~\ref{thm:always_sym} says that a sufficient condition for the MEP solutions to correspond to solutions to the original NEPv is that $g_i^Ty_{i+1} \neq  0$ for all $i$. The theory is inconclusive when $g_i^Ty_{i+1} = 0$ for some $i$. As we will illustrate now, this may lead to an MEP with a continuum of solutions, from which the NEPv solutions are difficult to identify. This would make our approach impractical for those choices.

Consider the NEPv \eqref{eq:nepv_general} with $m = 2$ and assume $g_1 = g_2$. We illustrate how choosing $g_1 = g_2$ leads to an MEP with a continuum of solutions. To this end, we first consider the following GEP parameterized with $\lambda$, with an additional orthogonality constraint
\begin{subequations}\label{eq:ortho_con_gep}
\begin{eqnarray}
-(A^T+\lambda B^T)y &=& (\mu_1 C_1^T+\mu_2 C_2^T)y,\\
g^Ty &=& 0.
\end{eqnarray}
\end{subequations}

Any solution $(\mu_1,\mu_2,y)$ to \eqref{eq:ortho_con_gep} leads to a corresponding symmetric solution $(\lambda,\mu_1,\mu_2,y,y,y)$ to the MEP given by
\begin{subequations}\label{eq:gsame_mep_t}
\begin{eqnarray}
-A^Ty_1 &=& (\lambda B^T+\mu_1 C_1^T+\mu_2 C_2^T)y_1\\
-(A^T+r_1g^T)y_2 &=& (\lambda B^T+\mu_1(C_1^T-s_1g^T)+\mu_2 C_2^T)y_2\\
-(A^T+r_2g^T)y_3 &=& (\lambda B^T+\mu_1 C_1^T+\mu_2 (C_2^T-s_2g^T))y_3.
\end{eqnarray}
\end{subequations}
Note that if $r_i \neq \alpha_i s_i$ for any $\alpha_i \in \mathbb{C}$, then every symmetric solution $(\gamma,\mu_1,\mu_2,y,y,y)$ to \eqref{eq:gsame_mep_t} also leads to a solution $(\mu_1,\mu_2,y)$ to \eqref{eq:ortho_con_gep} with $\lambda = \gamma$.
The MEP in \eqref{eq:gsame_mep_t} is the problem with left eigenvectors of \eqref{eq:nepv_general_linearizedt} with $m=2$ and $g_1 = g_2 = g$
\begin{subequations}\label{eq:gsame_mep}
\begin{eqnarray}
-Ax_1 &=& (\lambda B+\mu_1 C_1+\mu_2 C_2)x_1\\
-(A+gr_1^T)x_2 &=& (\lambda B+\mu_1(C_1-gs_1^T)+\mu_2 C_2)x_2\\
-(A+gr_2^T)x_3 &=& (\lambda B+\mu_1 C_1+\mu_2 (C_2-gs_2^T))x_3
\end{eqnarray}
\end{subequations}
Since \eqref{eq:ortho_con_gep} has more variables than equations, it can in general be solved for any $\lambda$ and hence, we have a continuum of solutions $(\lambda,\mu_1(\lambda),\mu_2(\lambda),y(\lambda),y(\lambda),y(\lambda))$ to \eqref{eq:gsame_mep_t}. Consequently, we have a continuum of solutions $(\lambda,\mu_1(\lambda),\mu_2(\lambda),x_1(\lambda),x_2(\lambda),x_3(\lambda))$ to \eqref{eq:gsame_mep}, which is undesirable. Hence, for the case $m = 2$, choosing equal $g_1$ and $g_2$ is problematic, and should be avoided.

\section{Numerical methods} \label{sec:num}
In the previous section, we showed how solutions of the nonlinear problem are contained in a linear multiparameter eigenvalue problem.
The eigenpairs can be extracted by solving a set of coupled generalized eigenvalue problems, but the large size of these problems is a limiting factor.
In this section we present two specialized methods that extract symmetric solutions by exploiting the structure of the problem.

\subsection{Residual inverse iteration}
The method of residual inverse iteration has been proposed for nonlinear eigenvalue problems in \cite{Neumaier:1985:RESINV} and was specialized for nonlinear two-parameter eigenvalue problems in \cite{plestenjak2016numerical}.
We use the notation from \cite{plestenjak2016numerical} to introduce the main properties of residual inverse iteration for the MEP defined by \eqref{eq:nepv_general_linearized}.
Furthermore, we show how we can exploit the structure of our problem to reduce computational cost and guarantee local convergence to a symmetric solution.
To ease the presentation, we first introduce the method for two-parameter problems and later generalize the results to problems with an arbitrary number of parameters.

First, we rewrite the 2EP \eqref{eq:nepv_general_linearized} using matrix-valued functions as
\begin{subequations}\label{eq:nepv_general_plestenjak}
    \begin{eqnarray}
    (A+\lambda B + \mu_1 C_1) x_1 &:=& T_1(\lambda,\mu_1)x_1=0 \\
    \left(A+g_1r_1^T + \lambda B + \mu_1(C_1-g_1s_1^T)\right)x_2 &:=& T_2(\lambda,\mu_1)x_2 = 0.
    \end{eqnarray}
\end{subequations}
The solutions of this system can be expressed as zeros of
\begin{equation} \label{eq:big_F}
F(x_1,x_2,\lambda,\mu) :=
\begin{bmatrix}
T_1(\lambda,\mu)x_1 \\
T_2(\lambda,\mu)x_2 \\
v_1^Tx_1 - 1 \\
v_2^Tx_2 - 1
\end{bmatrix}
\end{equation}
where the last two rows are normalization constraints.
The vectors $v_1$ and $v_2$ are chosen not orthogonal to an eigenvector.
Residual inverse iteration is derived from the quasi-Newton method applied to \eqref{eq:big_F}:
\begin{equation}\label{eq:quasi_newton}
\tilde{J}_F(x_1^{(k)},x_2^{(k)},\lambda^{(k)},\mu^{(k)})
\begin{bmatrix}
\Delta x_1 \\
\Delta x_2 \\
\Delta \lambda \\
\Delta \mu
\end{bmatrix}
= -
\begin{bmatrix}
T_1(\lambda^{(k)}, \mu^{(k)})x_1^{(k)}\\
T_2(\lambda^{(k)}, \mu^{(k)})x_2^{(k)}\\
v_1^Tx_1^{(k)} - 1 \\
v_2^Tx_2^{(k)} - 1
\end{bmatrix}
\end{equation}
with approximated Jacobian $\tilde{J}_F(x_1^{(k)},x_2^{(k)},\lambda^{(k)},\mu^{(k)})$ equal to
\begin{equation}\label{eq:resinv_jac}
\begin{bmatrix}
T_1(\sigma, \tau) & 0 & \frac{\partial T_1}{\partial \lambda}x_1^{(k)} & \frac{\partial T_1}{\partial \mu}x_1^{(k)} \\
0 & T_2(\sigma, \tau) & \frac{\partial T_2}{\partial \lambda}x_2^{(k)} & \frac{\partial T_2}{\partial \mu}x_2^{(k)} \\
v_1^T & 0 & 0 & 0 \\
0 & v_2^T & 0 & 0\\
\end{bmatrix}.
\end{equation}
The matrices $T_1(\sigma,\tau)$ and $T_2(\sigma,\tau)$ are nonsingular approximations to $T_1(\lambda^{(k)},\mu^{(k)})$ and $T_2(\lambda^{(k)},\mu^{(k)})$ respectively.
The arguments of the partial derivatives are omitted because they are constant for linear functions.
\begin{theorem}\label{thrm:resinv_equivalence}
    If $T_1$ and $T_2$ are linear and $x_1^{(k)}$ and $x_2^{(k)}$ are normalized in every iteration, the quasi-Newton system \eqref{eq:quasi_newton}-\eqref{eq:resinv_jac} is equivalent to the following update steps:
    \begin{equation}\label{eq:resinv_2x2}
    \begin{bmatrix}
    v_1^TT_1(\sigma,\tau)^{-1}\frac{\partial T_1}{\partial \lambda}x_1^{(k)} &
    v_1^TT_1(\sigma,\tau)^{-1}\frac{\partial T_1}{\partial \mu}x_1^{(k)} \\
    v_2^TT_2(\sigma,\tau)^{-1}\frac{\partial T_2}{\partial \lambda}x_2^{(k)} &
    v_2^TT_2(\sigma,\tau)^{-1}\frac{\partial T_2}{\partial \mu}x_2^{(k)}
    \end{bmatrix}
    \begin{bmatrix}
    \Delta \lambda \\
    \Delta \mu
    \end{bmatrix}
    = -\begin{bmatrix}
    \gamma_1\\
    \gamma_2
    \end{bmatrix}
    \end{equation}
    with $\gamma_i = v_i^TT_i(\sigma,\tau)^{-1}T_i(\lambda^{(k)},\mu^{(k)})x_i^{(k)}$ for $i=1,2$ and
    \begin{eqnarray} \label{eq:resinv_update}
    x_1 = x_1^{(k)} - T_1(\sigma,\tau)^{-1}T_1(\lambda^{(k+1)}, \mu^{(k+1)})x_1^{(k)}\\
    \label{eq:resinv_update2}
    x_2 = x_2^{(k)} - T_2(\sigma,\tau)^{-1}T_2(\lambda^{(k+1)}, \mu^{(k+1)})x_2^{(k)}
    \end{eqnarray}
    such that $x_1^{(k+1)} = \frac{x_1}{v_1^Tx_1}$ and $x_2^{(k+1)} = \frac{x_2}{v_2^Tx_2}$.
\end{theorem}
\begin{proof}
    Solving the first two rows of \eqref{eq:quasi_newton} for $\Delta x_i$ results in
    \begin{equation}\label{eq:Dxi}
        \Delta x_i = T_i(\sigma,\tau)^{-1}\left( T_i(\lambda^{(k)},\mu^{(k)})
        +\Delta\lambda \frac{\partial T_i}{\partial\lambda} + \Delta\mu\frac{\partial T_i}{\partial \mu} \right)x_i^{(k)} \quad i=1,2.
    \end{equation}
    If the vectors $x_1^{(k)}$ and $x_2^{(k)}$ are normalized in each iterations, it follows from the last two rows of \eqref{eq:quasi_newton} that $v_1^T\Delta x_1=0$ and $v_2^T\Delta x_2=0$.
    Solving $v_1^T\Delta x_1=0$ and $v_2^T\Delta x_2=0$ for $\Delta\lambda$ and $\Delta\mu$ leads to the $2\times2$ system \eqref{eq:resinv_2x2}.
    Notice that the matrix of this system is the Schur complement of \eqref{eq:resinv_jac} with respect to the bottom right block of zeros.
    Now, since $T_i$ is linear, it follows that
    \begin{equation}
    \Delta x_i = T_i(\sigma,\tau)^{-1}T_i(\lambda^{(k+1)},\mu^{(k+1)})x_i^{(k)} \quad i=1,2,
    \end{equation}
    which is equivalent to \eqref{eq:resinv_update}-\eqref{eq:resinv_update2}.
    A final normalization step, $x_i^{(k+1)} = \frac{x_i^{(k+1)}}{v_i^Tx_i^{(k+1)}}$ for $i=1,2$, proofs the stated equivalence.
\end{proof}

\begin{theorem}
    (Theorem 5.2 \cite{plestenjak2016numerical}) Let $(\lambda^{(*)},\mu^{(*)})$ be a simple eigenvalue of \eqref{eq:nepv_general_plestenjak} and let $x_1\otimes x_2$ be a corresponding right eigenvector such that $v_1^Tx_1 = v_2^Tx_2=1$. If $(\lambda^{(0)}, \mu^{(0)})$ is sufficiently close to the eigenvalue $(\lambda^{(*)},\mu^{(*)})$, then residual inverse iteration has linear convergence close to the solution.
\end{theorem}

\begin{theorem} \label{thm:symmetry_resinv}
    If residual inverse iteration is applied to problem \eqref{eq:nepv_general_plestenjak} with $v_1=v_2=v$ and $x_1^{(k)}=x_2^{(k)}=x^{(k)}$ such that $v^TT_1(\sigma,\tau)^{-1}g_1\neq 0$, then $x_1^{(k+1)}=x_2^{(k+1)}=x^{(k+1)}$ for $k=0,1,2,\dots$.
\end{theorem}
\begin{proof} \label{thrm:resinv}
    See proof of Theorem~\ref{thm:gen_symmetry_resinv}
\end{proof}

The result presented in Theorem~\ref{thm:symmetry_resinv} can be used to reduce the computational cost of the residual inverse iteration method for solving NEPv problems.
The structure exploiting method for finding symmetric solutions of \eqref{eq:nepv_general_plestenjak} is given in Algorithm~\ref{alg:plestenjak_symm}.
First, notice that in the linear case $\frac{\partial T_i}{\partial \lambda}$ and $\frac{\partial T_i}{\partial \mu}$ are constants such that the values $u_i$ and $w_i$ can be precalculated before the loop to reduce the computational cost.
Second, the symmetric method uses only one solution vector, $x$, which implies that we only need to update one vector.

\begin{algorithm}
    \caption{Symmetric residual inverse iteration for NEPv with 1 nonlinear term}
    \label{alg:plestenjak_symm}
    Start with $\sigma =\lambda^{(0)}, \tau=\mu^{(0)}, v$ and $x^{(0)}$  such that $v^Tx^{(0)}=1$.\\
    Compute $u_i=v^TT_i(\sigma,\tau)^{-1}\frac{\partial T_i}{\partial \lambda}$ for $i=1,2$.\\
    Compute $w_i=v^TT_i(\sigma,\tau)^{-1}\frac{\partial T_i}{\partial \mu}$ for $i=1,2$.\\
    \For{k=0,1,2,\dots} {
        Compute $\gamma_i = v^TT_i(\sigma,\tau)^{-1}T_i(\lambda^{(k)},\mu^{(k)})x^{(k)}$ for $i=1,2$.\\
        Solve
        \[
        \begin{bmatrix}
        u_1^Tx^{(k)} & w_1^Tx^{(k)} \\
        u_2^Tx^{(k)} & w_2^Tx^{(k)}
        \end{bmatrix}
        \begin{bmatrix}
        \Delta\lambda\\
        \Delta\mu
        \end{bmatrix}
        = -
        \begin{bmatrix}
        \gamma_1\\
        \gamma_2
        \end{bmatrix}.
        \]\\
        Update $\lambda^{(k+1)}=\lambda^{(k)}+\Delta\lambda$ and $\mu^{(k+1)}=\mu^{(k)}+\Delta\mu$.\\
        $x^{(k+1)} = x^{(k)} - T_1(\sigma,\tau)^{-1}T_1(\lambda^{(k+1)}, \mu^{(k+1)})x^{(k)}$\\
        Normalize $x^{(k+1)}=x^{(k+1)}/(v^Tx^{(k+1)})$.\\
    }
\end{algorithm}

\subsubsection*{Generalization to more parameters}
\newcommand{\mmu}{\boldsymbol{\mu}}
\newcommand{\ttau}{\boldsymbol{\tau}}
The theoretical derivation in the preceding part is deliberately limited to two parameters to keep the notation clear.
However, the method can be elegantly generalized to handle more than two parameters.
The derivation starts by rewriting the linearization \eqref{eq:nepv_general_linearized} as an MEP with $m+1$ parameters in the following form

\begin{subequations}\label{eq:gen_nepv_general_plestenjak}
    \begin{eqnarray}
    \left(A+\lambda B +\sum_{i=1}^{m} \mu_i C_i\right) x_1 &:=& T_1(\lambda,\mmu)x_1=0 \\
    \left(A+g_1r_1^T + \lambda B + \sum_{i=1}^{m} \mu_i C_i - \mu_1g_1s_1^T\right)x_2 &:=& T_2(\lambda,\mmu)x_2 = 0 \\
    \vdots \nonumber\\
    \left(A+g_mr_m^T + \lambda B +\sum_{i=1}^{m} \mu_i C_i  -\mu_mg_ms_m^T)\right)x_{m+1} &:=& T_{m+1}(\lambda,\mmu)x_{m+1}=0.
    \end{eqnarray}
\end{subequations}
with $\mmu = (\mu_1,\dots,\mu_m)$.
The solutions to this system can be expressed as zeros of
\begin{equation} \label{eq:gen_big_F}
F(x_1,\dots,x_{m+1},\lambda,\mmu) :=
\begin{bmatrix}
T_1(\lambda,\mmu)x_1 \\
\vdots \\
T_{m+1}(\lambda,\mmu)x_{m+1} \\
v_1^Tx_1 - 1 \\
\vdots \\
v_{m+1}^Tx_{m+1} - 1
\end{bmatrix}.
\end{equation}
The quasi-Newton method is now applied to \eqref{eq:gen_big_F}:
\begin{equation}\label{eq:gen_quasi_newton}
\tilde{J}_F(x_1^{(k)},\dots, x_{m+1}^{(k)},\lambda^{(k)},\mmu^{(k)})
\begin{bmatrix}
\Delta x_1 \\
\vdots \\
\Delta x_{m+1} \\
\Delta \lambda \\
\Delta \mmu
\end{bmatrix}
= -F(x_1^{(k)},\dots, x_{m+1}^{(k)},\lambda^{(k)},\mmu^{(k)})
\end{equation}
with approximated Jacobian $\tilde{J}_F(x_1^{(k)},\dots, x_{m+1}^{(k)},\lambda^{(k)},\mmu^{(k)})$ equal to the block matrix
\begin{equation}
\begin{bmatrix}
J_{1,1} & J_{1,2} \\
J_{2,1} & 0
\end{bmatrix}
\end{equation}
where the blocks are defined as
\begin{equation}
J_{1,1} =
\begin{bmatrix}
T_1(\sigma, \ttau) & & \\
& \ddots & \\
& & T_{m+1}(\sigma,\ttau)
\end{bmatrix},
\end{equation}
\begin{equation}
J_{1,2} =
\begin{bmatrix}
\frac{\partial T_1}{\partial\lambda}x_1^{(k)} & \dots & \frac{\partial T_1}{\partial\mu_m}x_1^{(k)} \\
\vdots & \ddots &\vdots \\
\frac{\partial T_{m+1}}{\partial\lambda}x_{m+1}^{(k)} & \dots & \frac{\partial T_{m+1}}{\partial\mu_m}x_{m+1}^{(k)}
\end{bmatrix} \,\text{and}
\end{equation}
\begin{equation}
J_{2,1}  =
\begin{bmatrix}
v_1^T & & \\
& \ddots & \\
& & v_{m+1}^T
\end{bmatrix}.
\end{equation}
The matrices $T_i(\sigma,\ttau)$ are fixed nonsingular approximations to $T_i(\lambda^{(k)},\mmu^{(k)})$.
\begin{theorem}
    If $T_i$ is linear and $x_i^{(k)}$ is normalized in every iteration for all $i=1,\dots,m+1$, then the quasi-Newton system \eqref{eq:gen_quasi_newton} is equivalent to the following update steps:
    \begin{equation}\label{eq:gen_resinv_system}
    \begin{bmatrix}
    \nu_{\lambda,1} & \nu_{\mu_1,1} & \dots & \nu_{\mu_m,1} \\
    \nu_{\lambda,2} & \nu_{\mu_1,2} & \dots & \nu_{\mu_m,2} \\
    \vdots & \vdots & \ddots & \vdots\\
    \nu_{\lambda,m+1} & \nu_{\mu_1,m+1} & \dots & \nu_{\mu_m,m+1}
    \end{bmatrix}
    \begin{bmatrix}
    \Delta \lambda \\
    \Delta \mu_1 \\
    \vdots \\
    \Delta \mu_m
    \end{bmatrix}
    = -
    \begin{bmatrix}
    \gamma_1 \\
    \gamma_2 \\
    \vdots \\
    \gamma_{m+1}
    \end{bmatrix}
    \end{equation}
    with $\nu_{\phi,i} = v_i^TT_i(\sigma,\ttau)^{-1}\frac{\partial T_i}{\partial \phi}x_i^{(k)}$ and $\gamma_i = v_i^TT_i(\sigma,\ttau)^{-1}T_i(\lambda^{(k)},\mmu^{(k)})x_i^{(k)}$ for $i=1,\dots,m+1$ and
    \begin{equation} \label{eq:gen_resinv_update}
    z_i = x_i^{(k)} - T_i(\sigma,\ttau)^{-1}T_i(\lambda^{(k+1)}, \mmu^{(k+1)})x_i^{(k)}
    \end{equation}
    such that $x_i^{(k+1)} = \frac{z_i}{v_i^Tz_i}$.
\end{theorem}

The proof is analogous to the proof of Theorem~\ref{thrm:resinv_equivalence} and is omitted for brevity.
We will now restate an important property of residual inverse iteration when applied to problem \eqref{eq:gen_nepv_general_plestenjak}.

\begin{theorem} \label{thm:gen_symmetry_resinv}
    If residual inverse iteration is applied to problem \eqref{eq:gen_nepv_general_plestenjak} with $v_1=\dots=v_{m+1}=v$ and $x_1^{(k)}=\dots= x_{m+1}^{(k)}=x^{(k)}$ such that $v^TT_1(\sigma,\ttau)^{-1}g_i\neq 0$ for $i=1,\dots,m$, then $x_1^{(k+1)}=\dots=x_{m+1}^{(k+1)}=x^{(k+1)}$ for $k=0,1,2,\dots$.
\end{theorem}
\begin{proof}
    The following proof shows that $x_1^{(k+1)} = x_i^{(k+1)}$ if $x_1^{(k)} = x_i^{(k)}$ for $i=2,\dots,m+1$.

    The bottom rows of \eqref{eq:gen_quasi_newton} imply that $v^T\Delta x_1=v^T\left(x_1^{(k+1)}-x_1^{(k)}\right) =0$ and $v^T\Delta x_i=v^T\left(x_i^{(k+1)}-x_i^{(k)}\right) =0$ if $x_1^{(k)}$ and $x_i^{(k)}$ are normalized.
    Assuming $x_1^{(k)}=x_i^{(k)}$, implies that
    \begin{equation}\label{eq:gen_v_equality}
    v^Tx_1^{(k+1)} = v^Tx_i^{(k+1)}.
    \end{equation}
    Now rewrite the first and $i$-th row of \eqref{eq:gen_quasi_newton} as
    \begin{align}\label{eq:gen_first_row}
    T_1(\sigma,\ttau)x_1^{(k+1)} =& -\left(\Delta\lambda B + \Delta\mu_1 C_1 + \cdots + \Delta\mu_m C_m - T_1(\sigma,\ttau) + T_1(\lambda^{(k)}, \mmu^{(k)})\right)x_1^{(k)} \\ \label{eq:gen_second_row}
    T_i(\sigma,\ttau)x_i^{(k+1)} =& -\Big(\Delta\lambda B + \Delta\mu_1C_1 +\cdots +\Delta\mu_{i-1} \left(C_{i-1}-g_{i-1}s_{i-1}^T\right) +\cdots  \\
    &+\Delta\mu_m C_m- T_i(\sigma,\ttau) + T_i(\lambda^{(k)}, \mmu^{(k)})\Big)x_i^{(k)}.\nonumber
    \end{align}
    Subtracting \eqref{eq:gen_first_row} from \eqref{eq:gen_second_row} gives
    \begin{equation} \label{eq:gen_zeta_eq}
    T_1(\sigma,\ttau)(x_i^{(k+1)} - x_1^{(k+1)}) = g_{i-1}\zeta
    \end{equation}
    with
    \begin{equation}
    \begin{split}
    \zeta=s_{i-1}^Tx_{i-1}^{(k)}\Delta\mu_{i-1} + (r_{i-1}-\tau_{i-1} s_{i-1})^Tx_{i-1}^{(k)} - (r_{i-1}-\mu^{(k)}s_{i-1})^Tx_{i-1}^{(k)} - \\ (r_{i-1}-\tau_{i-1} s_{i-1})^Tx_2^{(k+1)}.
    \end{split}
    \end{equation}
    Finally, we left-multiply with $v^TT_1(\sigma,\ttau)^{-1}$:
    \begin{equation}
    v^T\left(x_i^{(k+1)} - x_1^{(k+1)}\right) = v^TT_1(\sigma,\ttau)^{-1}g_{i-1}\zeta.
    \end{equation}
    From \eqref{eq:gen_v_equality} we know that the left-hand side is zero.
    The assumption $v^TT_1(\sigma,\tau)^{-1}g_{i-1} \neq 0$ implies that $\zeta=0$.
    Since $T_1(\sigma,\ttau)$ is nonsingular, \eqref{eq:gen_zeta_eq} now leads to the conclusion that $x_1^{(k+1)} = x_i^{(k+1)}$.
\end{proof}

The result presented in Theorem~\ref{thm:gen_symmetry_resinv} can again be used to reduce the computational cost of the residual inverse iteration method for solving NEPv problems.
The structure exploiting method for finding symmetric solutions of \eqref{eq:gen_nepv_general_plestenjak} is given in Algorithm~\ref{alg:gen_plestenjak_symm}.

\begin{algorithm}
    \caption{Symmetric residual inverse iteration for NEPv with $m$ nonlinear terms}
    \label{alg:gen_plestenjak_symm}
    Start with $\sigma =\lambda^{(0)}, \ttau=\mmu^{(0)}, v$ and $x^{(0)}$ such that $v^Tx^{(0)}=1$.\\
    Compute all $\psi_{\phi,i}^T = v^TT_i(\sigma,\ttau)^{-1}\frac{\partial T_i}{\partial \phi}$ for $\phi\in(\lambda,\mu_1,\dots,\mu_m)$ and $ i\in(1,\dots,m+1)$.\\
    \For{k=0,1,2,\dots} {
        Compute $\gamma_i = v_i^TT_i(\sigma,\ttau)^{-1}T_i(\lambda^{(k)},\mmu^{(k)})x_i^{(k)}$ for $i\in(1,\dots,m+1)$.\\
        Solve
        \begin{equation*}
        \begin{bmatrix}
        \psi_{\lambda,1}^Tx^{(k)} & \psi_{\mu_1,1}^Tx^{(k)} & \dots & \psi_{\mu_m,1}^Tx^{(k)} \\
        \psi_{\lambda,2}^Tx^{(k)} & \psi_{\mu_1,2}^Tx^{(k)} & \dots & \psi_{\mu_m,2}^Tx^{(k)} \\
        \vdots & \vdots & \ddots & \vdots\\
        \psi_{\lambda,m+1}^Tx^{(k)} & \psi_{\mu_1,m+1}^Tx^{(k)} & \dots & \psi_{\mu_m,m+1}^Tx^{(k)}
        \end{bmatrix}
        \begin{bmatrix}
        \Delta \lambda \\
        \Delta \mu_1 \\
        \vdots \\
        \Delta \mu_m
        \end{bmatrix}
        =
        -
        \begin{bmatrix}
        \gamma_1 \\
        \gamma_2 \\
        \vdots \\
        \gamma_{m+1}
        \end{bmatrix}.
        \end{equation*}\\
        Update $\lambda^{(k+1)}=\lambda^{(k)}+\Delta\lambda$ and $\mu_i^{(k+1)}=\mu_i^{(k)}+\Delta\mu_i$ for $i\in(1,\dots,m+1)$.\\
        Set $x^{(k+1)} = x^{(k)} - T_1(\sigma,\ttau)^{-1}T_1(\lambda^{(k+1)}, \mmu^{(k+1)})x^{(k)}$.\\
        Normalize $x^{(k+1)}=x^{(k+1)}/(v^Tx^{(k+1)})$.\\
    }
\end{algorithm}

\subsection{Inverse iteration methods for operator determinants}
The method of inverse iteration for the standard eigenvalue problem is known to be globally convergent under very mild conditions, which is beneficial if no prior knowledge of the problem is available.
The application of this method to the matrix pencil of operator determinants $(\Delta_1, \Delta_0)$ allows for the extraction of nonlinear eigenpairs $(\lambda,x)$.
In the remainder of this section we discuss how this method can preserve symmetry in the eigenvector estimates while improving the convergence speed.

The goal is to find eigenpairs $(\lambda, z)$ of the generalized eigenvalue problem
\begin{equation}
\Delta_1 z = \lambda \Delta_0z
\end{equation}
with $\Delta_0,\Delta_1\in \mathbb{C}^{n^{m+1}\times n^{m+1}}$ the operator determinants derived from the NEPv \eqref{eq:nepv_general} with $m$ nonlinear terms and $z=x\otimes\cdots\otimes x\in\mathbb{C}^{n^{m+1}}$.
Unless stated otherwise, we assume that the operator determinants $\Delta_0$ and $\Delta_1$ are nonsingular.
In the general case, this problem has $n^{m+1}$ linearly independent eigenvectors, but, for the corresponding NEPv, only the $N_s=\begin{pmatrix}n+m\\m+1\end{pmatrix}$ symmetric eigenvectors are of interest.

A well known method to solve the eigenvalue problem is inverse iteration for which the iterates are dictated by
\begin{subequations}\label{eq:inv_iter}
\begin{eqnarray}
\left(\Delta_1-\sigma\Delta_0\right)\tilde{z}_{k+1} = \Delta_0z_k\\
z_{k+1} = \frac{\tilde{z}_{k+1}}{\lVert\tilde{z}_{k+1}\rVert} \label{eq:inv_iter_z}
\end{eqnarray}
\end{subequations}
with $\sigma\in\mathbb{C}$ an estimate of the desired eigenvalue such that $\Delta_1-\sigma\Delta_0$ is nonsingular.
Although this method will generally converge to the eigenvector corresponding to the eigenvalue closest to $\sigma$, it is possible to restrict convergence to symmetric eigenvectors only.

\begin{lemma}
    \label{lmm:symm_inv_s}
   If the pencil $(\Delta_1, \Delta_0)$ is regular and non-defective, $\sigma$ is not an eigenvalue and if $z_k$ is symmetric and decomposable as $z_k=\sum_{j=1}^{l}\alpha_j x_j\otimes\dots\otimes x_j$ for some index $l\leq N_s$ with $x_j$ an eigenvector of the NEPv, then $z_{k+1}$ in \eqref{eq:inv_iter_z} is symmetric too.
\end{lemma}
\begin{proof}
    If $z_k=\sum_{j=1}^{l}\alpha_j x_j\otimes\dots\otimes x_j$ then it follows from inverse iteration property on eigenvectors that $z_{k+1}=\sum_{j=1}^{l}\frac{\alpha_j}{\lambda_j - \sigma}x_j\otimes\dots\otimes x_j$ with $\lambda_j$ the eigenvalue related to the eigenvector $x_j$. It is then clear that $z_{k+1}$ is also  symmetric.
\end{proof}

This lemma leads to a more powerful result as described in the following theorem.
\begin{theorem} \label{thrm:symm_inv}
    If the NEPv has $N_s$ distinct solutions, the pencil $(\Delta_1, \Delta_0)$ is regular and non-defective, $\sigma$ is not an eigenvalue and if $z_k$ is symmetric, then $z_{k+1}$ is symmetric.
\end{theorem}
\begin{proof}
    If the NEPv has $N_s$ solutions, Lemma \ref{lmm:symmetric_vectors} dictates that if the pencil $(\Delta_1, \Delta_0)$ is regular and non-defective, there are $N_s$ linearly independent eigenvectors spaning the complete space of symmetric vectors that is of dimension $N_s$.
    This means that every symmetric vector $z_k$ can be decomposed as $z_k=\sum_{j=1}^{l}\alpha_j x_j\otimes\dots\otimes x_j$ for some index $l\leq N_s$ and with $x_j$ an eigenvector of the NEPv.
    By Lemma~\ref{lmm:symm_inv_s} the vector $z_{k+1}$ is then symmetric too, proving the stated result.
\end{proof}

The previous theorem not only shows how symmetric eigenvectors can be extracted using inverse iteration, but also suggests how the convergence rate only depends on those eigenvectors.
It is well known that, for general matrices, the convergence rate of inverse iteration is dictated by $\mathcal{O}\left(\left\vert\frac{\sigma - \lambda_{1*}}{\sigma - \lambda_{2*}}\right\vert^k \right)$ with $\lambda_{1*}$ and $\lambda_{1*}$ respectively the closest and second closest eigenvalue to $\sigma$.
The following theorem shows how this rate can be specified for symmetric starting vectors.
\begin{theorem}
    If the NEPv has $N_s$ distinct solutions, the pencil $(\Delta_1, \Delta_0)$ is regular and non-defective, $\sigma$ is not an eigenvalue and if $z_0$ is symmetric, then the convergence rate of inverse iteration \eqref{eq:inv_iter} is dictated by $\mathcal{O}\left(\left\vert\frac{\sigma - \lambda_1}{\sigma - \lambda_2}\right\vert^k \right)$ with $\lambda_1$ and $\lambda_2$ respectively the closest and second closest symmetric eigenvalue to $\sigma$, assuming that $z_0$ has a nonzero component in the corresponding eigenvectors.
\end{theorem}
\begin{proof}
    It follows from Lemma~\ref{lmm:symm_inv_s} and Theorem~\ref{thrm:symm_inv} that $z_{k+1}=\sum_{j=1}^{N_s}\frac{\alpha_j}{\lambda_j - \sigma}x_j\otimes\dots\otimes x_j$ with the index $j$ defined such that $\vert\lambda_1-\sigma\vert\leq\vert\lambda_2-\sigma\vert\leq\cdots\leq\vert\lambda_{N_s}-\sigma\vert$.
    This means that after each iteration the component of the first eigenvector grows by a factor $\frac{1}{\vert\sigma-\lambda_1\vert}$, while the component of the second only grows by a factor $\frac{1}{\vert\sigma-\lambda_2\vert}$.
    The coefficient of $x_1\otimes\cdots\otimes x_1$ will therefore grow faster than the coefficient of the eigenvector of the second closest symmetric eigenvalue by a factor of $\left\vert\frac{\sigma - \lambda_2}{\sigma - \lambda_1}\right\vert$ which leads to a convergence rate of $\mathcal{O}\left(\left\vert\frac{\sigma - \lambda_1}{\sigma - \lambda_2}\right\vert^k \right)$.
\end{proof}

As a result of this theorem, the convergence of inverse iteration is independent of the nonsymmetric solutions of the linearization if the starting vector is symmetric.

This section concludes with a discussion about the complexity of inverse iteration.
If the nonlinear problem of dimension $n$ has $m$ nonlinear terms, each iteration requires solving the linear system  given in \eqref{eq:inv_iter} that is of dimension $n^{m+1}\times n^{m+1}$.
A standard implementation of inverse iteration would therefore lead to a computational complexity of $\mathcal{O}(n^{3(m+1)})$.
However, it is possible to reduce the computational cost by exploiting the structure of the operator determinants \cite{meerbergen2015sylvester}.
To this end, we first introduce the definition of the $\text{vec}(\cdot)$ operator.
\begin{definition}
    The vectorization of a tensor is the bijective mapping that reorders the elements of the tensor $\mathcal{A}\in\mathbb{C}^{n_1\times\cdots\times n_d}$ in a vector such that
    \begin{equation*}
        {\normalfont \text{vec}}(\mathcal{A}) = [a_{1,1,\ldots,1} \enspace a_{2,1,\ldots,1} \enspace\cdots\enspace a_{n_1,1,\ldots,1} \enspace a_{1,2,1,\ldots,1} \enspace a_{2,2,1,\ldots,1} \enspace\cdots\enspace a_{n_1,n_2,\ldots,n_d}]^T.
    \end{equation*}
\end{definition}

Consider for simplicity the case where $m=1$ for which \eqref{eq:inv_iter} is rewritten as
\begin{equation} \label{eq:orig_explicit}
\begin{split}
\left(C\otimes (A+g_1r_1^T+\sigma B) - (A+\sigma B)\otimes(C-g_1s_1^T) \right)z_{k+1} \\
= \left( B\otimes(C-g_1s_1^T) - C\otimes B \right)z_k.
\end{split}
\end{equation}
This equation is transformed to a Sylvester equation by using the ``vec trick'' \cite{roth1934direct}:
\begin{equation}\label{eq:sylvester}
\begin{split}
(A+g_1r_1^T+\sigma B)Z_{k+1}C^T - (C-g_1s_1^T)Z_{k+1}(A+\sigma B)^T \\
= (C-g_1s_1^T)Z_kB^T-BZ_kC^T
\end{split}
\end{equation}
with $\text{vec}(Z_{k+1}) = z_{k+1}$ and $\text{vec}(Z_{k}) = z_{k}$.
The computational cost reduces from $\mathcal{O}(n^6)$ for the standard implementation to $\mathcal{O}(n^3)$ for the Sylvester equation approach.
This computational complexity is based on the Bartels-Stewart method \cite{Bartels:1972:LYAP} to solve the Sylvester equation.
For $m>1$, equation \eqref{eq:inv_iter} can be transformed to a tensor Sylvester equation by using the following generalization of the vec trick:
\begin{equation}
(T_1\otimes\dots\otimes T_{m+1})\text{vec}(\mathcal{Y}) = \text{vec}( \mathcal{Y}\times_1T_{m+1}\times_2\dots\times_{m+1}T_1 )
\end{equation}
where $\mathcal{Y}$ is a square tensor of dimension $n$ and order $m+1$, $T_i\in\mathbb{C}^{n\times n}$ for $i=1,\ldots,m+1$, and the operator $\times_i$ is the $i$-mode tensor product described by \cite{kolda2009tensor}.
The numerical methods that solve these tensor Sylvester equations \cite{ali2016krylov,chen2012projection,huang2020iterative,kressner2010krylov,ruymbeek2020subspace} are based on iterative schemes which makes it harder to quantify the computational cost reduction for the Sylvester approach for $m>1$.
A more thorough discussion of the existing methods falls outside of the scope of this paper.

\subsection{Numerical experiments}
This section shows the convergence properties of three iterative methods: residual inverse iteration (RI) as in \cite{plestenjak2016numerical}, residual inverse iteration with symmetry exploitation (RIS) as in Algorithm~\ref{alg:gen_plestenjak_symm} and inverse iteration (II) defined by \eqref{eq:inv_iter}.
These methods all require initial estimates of the desired eigenpair $(\lambda,x)$.
RI and RIS require additional estimates of the remaining components of the eigenvalue $(\mu_1,\dots,\mu_m)$ that are calculated using the eigenvector estimate $x^{(0)}$ as $\tau_i = f_i(x^{(0)})$.
Finally, we assume that the normalization vectors are $v=v_1=v_2=v_3$ and the starting vectors are $x^{(0)} = x_1^{(0)} = x_2^{(0)} = x_3^{(0)}$ for RI and RIS while the starting vector for II is $z_0=x^{(0)}\otimes x^{(0)}$ for Example~1 and $z_0=x^{(0)}\otimes x^{(0)}\otimes x^{(0)}$ for Example~2.

\subsubsection{Example 1}
\begin{figure}
    \centering
    \includegraphics[width=0.8\textwidth]{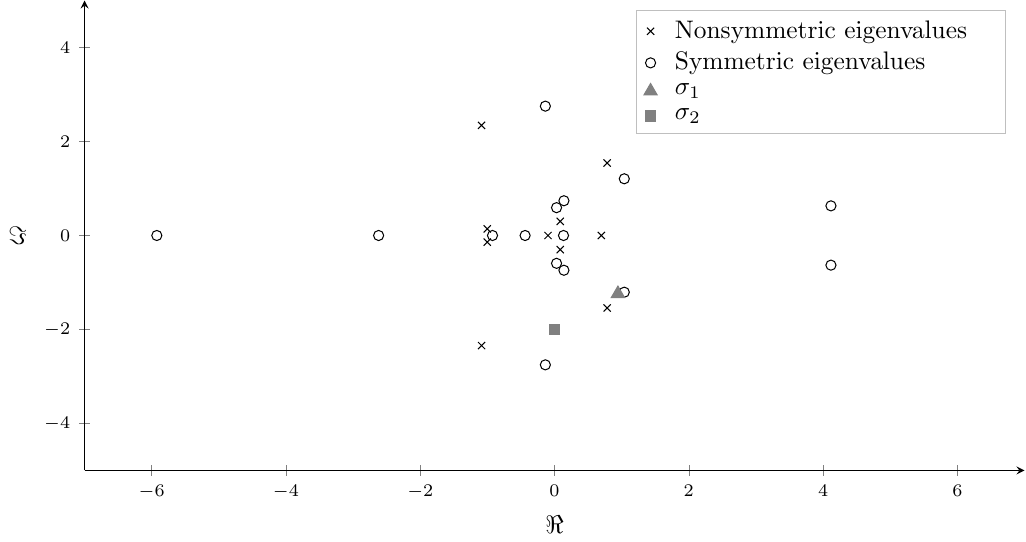}
    \caption{Eigenvalues of 2EP in Example 1.}
    \label{fig:ex1_eigs}
\end{figure}

Consider the following problem
\[
(A+\lambda B + f(x)C)x=0
\]
defined by randomly generated matrices $A,B,C \in \mathbb{R}^{5\times 5}$ and by the function $f(x)=\frac{r^Tx}{s^Tx}$ with $r,s \in \mathbb{R}^5$ also randomly generated.
The corresponding two parameter eigenvalue problem is given by
\[
\begin{cases}
-Ax_1 = (\lambda B+\mu C)x_1 \\
-(A+gr^T)x_2 = (\lambda B + \mu(C-gs^T))x_2
\end{cases}
\]
with $g\in\mathbb{R}^5$ randomly generated.
Due to the small problem size, it is feasible to solve the 2EP using the explicit operator determinant formulation $\Delta_1z=\lambda\Delta_0z$ for which the eigenvalues are given in Figure~\ref{fig:ex1_eigs}.
There are 15 eigenvalues related to symmetric eigenvectors and 10 eigenvalues related to nonsymmetric eigenvectors, which is to be expected for a $5\times5$ nonlinear problem.
Using the results from the operator determinant eigenvalues, we can now show the convergence properties of every method.

\tikzsetnextfilename{ex1_eigenvalues}

\begin{figure}
    \centering
    \tikzsetnextfilename{ex1_conv1}
    \begin{subfigure}[b]{0.48\textwidth}
        \includegraphics[width=\textwidth]{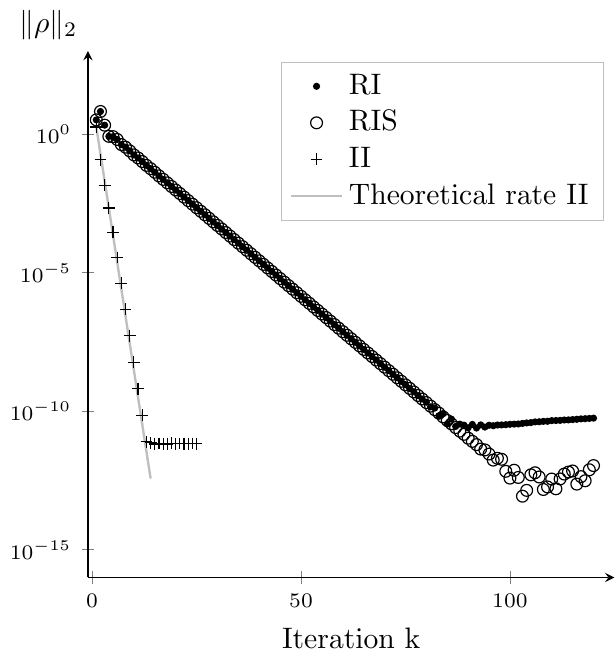}
        \caption{$\sigma=0.9440 - 1.2382\imath$}
        \label{fig:ex1_conv1}
    \end{subfigure}
    ~
    \begin{subfigure}[b]{0.48\textwidth}
        \tikzsetnextfilename{ex1_conv2}
        \includegraphics[width=\textwidth]{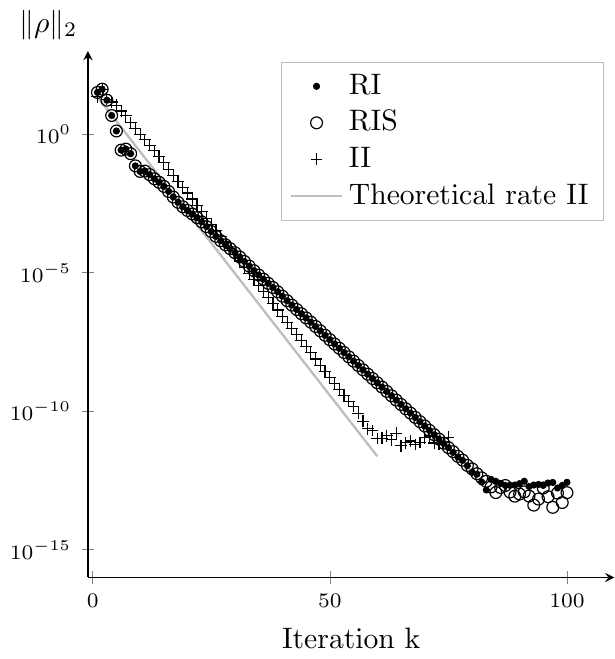}
        \caption{$\sigma=-2\imath$}
        \label{fig:ex1_conv2}
    \end{subfigure}
    \caption{Residual of NEPv of Example 1.}\label{fig:ex1_conv}
\end{figure}

The first test assumes a priori knowledge about the desired eigenpair: $\sigma=0.9440 - 1.2382\imath$ is close to a desired eigenvalue and $x^{(0)}$ is calculated as the corresponding eigenvector with a random perturbation of norm $10^{-1}$.
The starting value $\sigma$ is indicated in Figure~\ref{fig:ex1_eigs} as a triangle.
Figure~\ref{fig:ex1_conv1} shows the residual, $\rho$, of the NEPv and shows the convergence behaviour of the three methods.
It is clear that RI and RIS are mathematically equivalent by looking at the first part of the convergence, but the last part shows numerical rounding errors that break the symmetry in RI which causes it to diverge from RIS.
The convergence plot for II shows that it stagnates around $10^{-11}$, which can be explained by looking at the condition number of $\Delta_0$ and $\Delta_1$ which is $\mathcal{O}(10^4)$.
This means that \eqref{eq:inv_iter} cannot be solved with a smaller residual in double precision.
Lastly, note that II follows the symmetric convergence speed of $\mathcal{O}\left(0.11^k\right)$ instead of the nonsymmetric estimate $\mathcal{O}\left(0.29^k\right)$.

The second test for this example assumes no prior knowledge and aims to find the eigenvalue closest to $\sigma = -2\imath$, depicted as a square in \ref{fig:ex1_eigs}.
All three methods start with the same randomly generated vector $x^{(0)}$, which leads to the convergence result as given in Figure~\ref{fig:ex1_conv2}.
The method of II follows the predicted convergence rate and stagnates at the same level as in the first test.
RI and RIS converge linearly to a solution, but closer inspection shows that they converged to $\lambda\approx0.14 - 0.74\imath$ which is not the eigenvalue nearest to $\sigma$.
We observe linear convergence for all tested methods even if RI and RIS are only guaranteed to converge locally.

\subsubsection{Example 2}
Consider the following NEPv with two nonlinear terms
\begin{equation}
\begin{split}
\left( A + \lambda B + f_1(x)C_1 + f_2(x)C_2)\right)x=0 \\
f_1(x) = \frac{r_1^Tx}{s_1^Tx}, f_2(x) = \frac{r_2^Tx}{s_2^Tx}\\
\end{split}
\end{equation}
which is linearized as the following 3-parameter eigenvalue problem
\[
\begin{cases}
-Ax_1 = (\lambda B+\mu_1 C_1+\mu_2C_2)x_1 \\
-(A+g_1r_1^T)x_2 = (\lambda B + \mu_1(C_1-g_1s_1^T) + \mu_2C_2)x_2 \\
-(A+g_2r_2^T)x_3 = (\lambda B + \mu_1C_1 + \mu_2(C_2-g_2s_2^T))x_3.
\end{cases}
\]
In this example, the matrices $A,B,C_1,C_2\in\mathbb{R}^{10\times10}$ and vectors $r_1,r_2,s_1,s_2,g_1,g_2\in\mathbb{R}^{10}$ are all randomly generated.
The size of the NEPv is rather small such that the equivalent generalized eigenvalue problem is still feasible to solve.
This way, a generic eigenproblem solver can be used to find all symmetric eigenpairs and to compare them to the result from the iterative methods.
Similarly as in Example~1, we present two situations: the first will assume that prior knowledge of the eigenpair is available, the second situation will assume no knowledge of the solution.

\begin{figure}
    \centering
    \tikzsetnextfilename{ex2_conv1}
    \begin{subfigure}[b]{0.48\textwidth}
        \includegraphics[width=\textwidth]{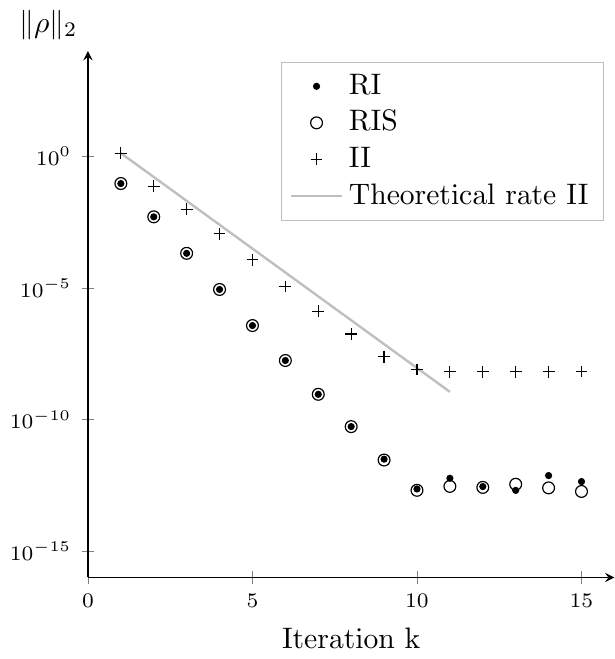}
        \caption{$\sigma=-1.2769 - 0.0442\imath$}
        \label{fig:ex2_conv1}
    \end{subfigure}
    ~
    \tikzsetnextfilename{ex2_conv2}
    \begin{subfigure}[b]{0.48\textwidth}
        \includegraphics[width=\textwidth]{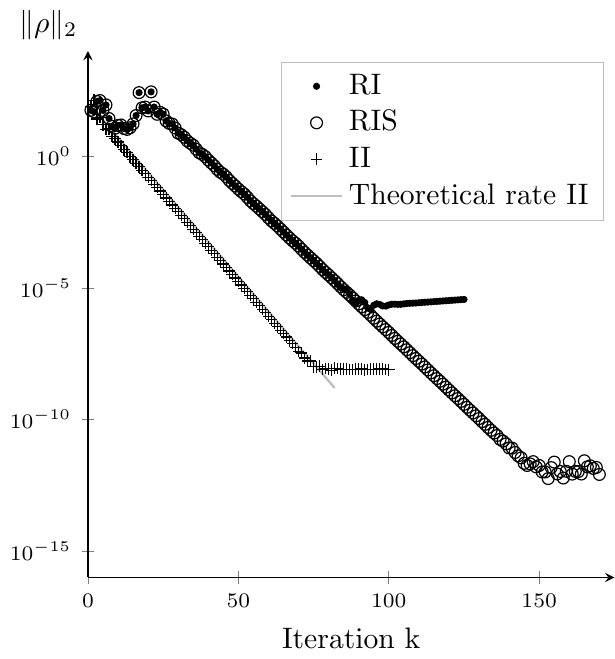}
        \caption{$\sigma=-4-3\imath$}
        \label{fig:ex2_conv2}
    \end{subfigure}
    \caption{Residual of NEPv of Example 2.}\label{fig:ex2_conv}
\end{figure}

First consider the case where prior knowledge of the problem is available to start the iteration.
We use the result from the operator determinant eigenvalue problem to pick an eigenvalue $\lambda \approx -1.2830 - 0.0363\imath$ and its corresponding normalized eigenvector $x$.
The starting values for the iterative methods are a perturbation of this eigenpair $\sigma = -1.2769 - 0.0442\imath$ and $x^{(0)} = x + e$ with $e$ a randomly generated vector of norm $10^{-2}$.
The starting value for the eigenvector is used to determine the values of $\mu_1^{(0)}=f_1(x^{(0)})$ and $\mu_2^{(0)}=f_2(x^{(0)})$ .
Figure~\ref{fig:ex2_conv1} shows the residual as a function of the number of iterations.
As expected, all three methods converge to the closest eigenpair within a small number of iterations.
The method of inverse iterations follows the symmetric convergence rate of $\mathcal{O}(0.12^k)$ instead of the standard rate $\mathcal{O}(0.22^k)$.
Again, note that II stagnates at the level around $10^{-9}$ which can be explained by the same reasoning as in Example 1 since the condition number of $\Delta_0$ and $\Delta_1$ is $\mathcal{O}(10^7)$.

The second part of this example considers the situation in which no prior knowledge is available and we initialize with a target $\sigma=-4-3\imath$ and a random starting vector while the initial values for $\mu_1^{(0)}$ and $\mu_2^{(0)}$ are calculated in the same way as in the previous part of this example.
The residual norms are plotted in Figure~\ref{fig:ex2_conv2}.
Again, RI and RIS show the same convergence behaviour but, between iterations 80 and 90, RI diverges from RIS again due to numerical rounding errors.
However, since RIS did not start sufficiently close to an eigenpair, it converged to the eigenvalue $\lambda\approx-2.6086 + 0.2827\imath$ and not to the closest eigenvalue $\lambda\approx-3.9910 - 2.2599\imath$.
The method of inverse iteration did converge, with limited accuracy,  to the closest eigenpair and followed the theoretical symmetric convergence rate of $\mathcal{O}(0.74^k)$ which is better than the standard convergence rate of $\mathcal{O}(0.94^k)$.

\subsubsection{Example 3}
Consider the following differential equation in one spatial dimension
\begin{equation}
   u^{\prime\prime}(x) +\lambda k_1(x)u(x) + f(u)k_2(x)u(x)=0
\end{equation}
defined for the domain $x\in[-1,1]$ with $u(-1)=u(1)=0$.
The function $f(u)=\frac{\alpha(u)}{\beta(u)}$ is defined by $\alpha(u) = \int_{-1}^{1} g(x)u(x) dx$ and $\beta(u) = u^{\prime}(0)$ with $g(x) = e^{-\gamma x^2}$.
For this example, we set $k_1(x) = 1+\frac{1}{2}\text{tanh}(5x)$, $k_2(x) = 1+\frac{1}{2}\text{cos}(\pi x)$ and $\gamma=10$.
The solutions of this differential equation can be interpreted as vibrations in a medium with a damping coefficient that depends on both the spatial location and the shape of the vibrations itself.
By spatial discretization, the equation can be written as a NEPv:
\begin{equation} \label{eq:disc_pde}
    Av + \lambda K_1v + \frac{a^Tv}{b^Tv}K_2v = 0
\end{equation}
with $v$ the discretization of $u(x)$ and where the derivatives are discretized by finite differences and the integral by the trapezoidal rule.
The NEPv \eqref{eq:disc_pde} is of the general form \eqref{eq:nepv_general}-\eqref{eq:nepv_general_f} and can be linearized to a MEP.

\begin{figure}
    \centering
    \tikzsetnextfilename{ex3_conv}
    \begin{subfigure}[b]{0.48\textwidth}
        \includegraphics[width=\textwidth]{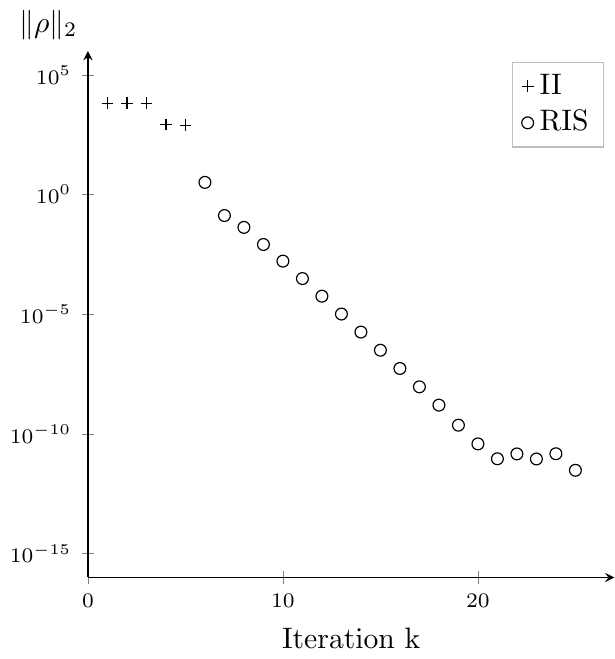}
        \caption{Residual}
        \label{fig:ex3_conv}
    \end{subfigure}
    ~
    \begin{subfigure}[b]{0.48\textwidth}
        \tikzsetnextfilename{ex3_result}
        \includegraphics[width=\textwidth]{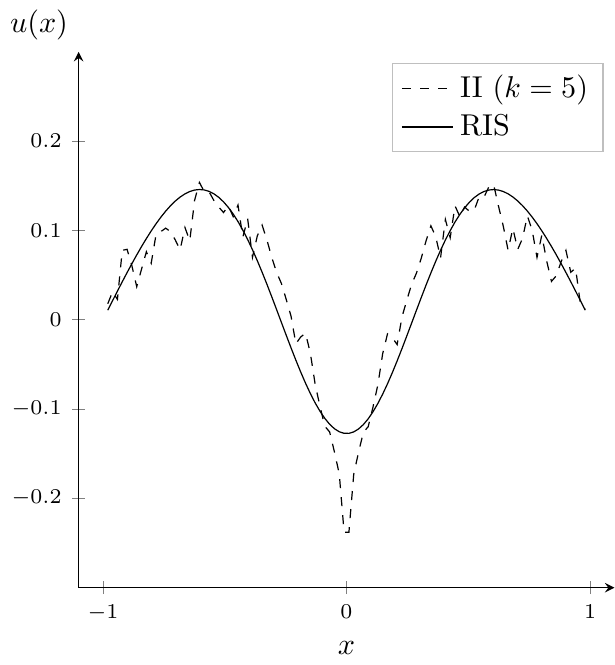}
        \caption{Calculated solution}
        \label{fig:ex3_result}
    \end{subfigure}
    \caption{Convergence and solution for Example~3.}\label{fig:ex3}
\end{figure}

Using $100$ discretization points, we will show how inverse iteration and residual inverse iteration can go hand in hand to find desired solutions in a reliable and fast way.
From the theory and the first two examples, we know that inverse iteration is a very robust, but more expensive, method while residual inverse iteration requires better starting values but is cheaper.
Therefore, we use a hybrid approach for this example: start with a small number of iterations of inverse iteration to improve the starting values and continue with residual inverse iteration until convergence.
The result is presented in Figure~\ref{fig:ex3} where we use $\sigma = 1$ and a random symmetric starting vector for II.
The intermediate approximation after five iterations of II is used as a starting value for RIS.
The approximation of $u(x)$ after II is shown together with the final result after RIS in Figure~\ref{fig:ex3_result}.
These results show how the proposed methods can be complementary in practical applications: II reliably finds a rough estimate of the desired solution, and RIS cheaply converges to the solution with small residual.

\section{Discussion} \label{sec:concl}
We introduced a new way of analyzing the NEPv by the concept of linearization for eigenvector nonlinearities.
The linearization is derived for a subset of nonlinear problems and introduces some additional parameters that can be chosen freely.
We show how to linearize problems of this class and characterize the number of solutions. We also derive conditions that guide us on how to choose the free parameters. Contrary to existing methods for NEPv, the proposed linearization is able to identify any and all solutions, under mild assumptions.

The new linearization technique also allows for the use and adaption of well-understood state-of-the-art methods.
We proposed two structure-exploiting numerical schemes: residual inverse iteration and inverse iteration.
Residual inverse iteration requires less calculations per iteration but needs accurate starting values to obtain predictable convergence.
Inverse iteration is more expensive but is in general globally convergent to the desired solution.
Three numerical examples confirmed these theoretical results and showed that the two methods are complementary: a small number of iterations of inverse iteration lead to a good starting value that can be used for residual inverse iteration.
A thorough discussion of the implementation of the numerical methods falls outside the scope of this paper, but it is an interesting point for future research. 

\section*{Acknowledgements}
We thank Bor Plestenjak (Univ. Ljubljana) for explanations
concerning general MEPs, Lars Karlsson (Univ. Ume\aa) for
comments concerning GEP pencil singularity, and Simon Telen (Max Planck Institute, Liepzig) for insights in the multi-homogenous B\'ezout theorem.

The work by Rob Claes and Karl Meerbergen is supported by the Research Foundation – Flanders (FWO) Grant G0B7818N and the KU Leuven Research Council.

\bibliographystyle{plain}
\bibliography{fulljabref}
\end{document}